\documentclass[a4paper,11pt]{article}

\usepackage{amsmath}
\usepackage{amsthm}
\usepackage{amssymb}

\usepackage[latin1]{inputenc}

\usepackage[dvips=true,%
a4paper=true,%
pdftitle={titre},%
pdfauthor={Jean-Francois Burnol},%
pdfstartview=FitH,%
pdfpagemode=None%
]{hyperref}

\def\ladate{July 19, 2010}

\newtheorem{theorem}{Theorem}

\newtheorem{proposition}[theorem]{Proposition}

\theoremstyle{definition}

\newtheorem{remark}{Remark}

 \newcommand{\RR}{{\mathbb R}}
 \newcommand{\CC}{{\mathbb C}}

\newcommand{\cE}{{\mathcal E}}
\newcommand{\cF}{{\mathcal F}}

\DeclareMathOperator{\Imag}{Im}

\begin{document}

\title{Hilbert spaces of entire functions with trivial zeros}

\author{Jean-Fran\c cois Burnol}

\date{\ladate}

\maketitle

\begin{center}
\begin{small}

\begin{abstract}   Let $H$ be a Hilbert space of entire functions. Let $H'$ be
the space of the functions $f(z)/\prod_i(z-z_i)$ where $f$ belongs to $H$ and
vanishes at $n$ given complex points $z_i$. We compute a suitable $E$ function
for $H'$  when one is given for $H$.
\end{abstract}

\bigskip

\begin{quote}
Universit\'e Lille 1\\ 
UFR de Math\'ematiques\\ 
Cit\'e scientifique M2\\ 
F-59655 Villeneuve d'Ascq\\ 
France\\
burnol@math.univ-lille1.fr\\
\end{quote}
\end{small}
\end{center}

\setlength{\normalbaselineskip}{16pt}
\baselineskip\normalbaselineskip
\setlength{\parskip}{6pt}

\section{Hilbert space of entire functions with vanishing conditions}

Let $H$ be an Hilbert space, whose vectors are entire functions, and such that
the evaluations at complex numbers are continuous linear forms, hence
correspond to specific vectors $Z_z$:
\begin{equation}\label{eq:ev}
Z_z\in H,\quad \forall f\in H\qquad (Z_z,f) = f(z)
\end{equation} Scalar products $(g,f)$ in this paper will be complex linear in
$f$  and conjugate linear in $g$. Using the Banach-Steinhaus theorem we know
that evaluations of derivatives  $f\mapsto f^{(k)}(z)$
are also bounded linear forms on the Hilbert space $H$. We will write $Z_{z,k}$ for the corresponding
vectors.

If $H$ satisfies the axiomatic framework of \cite{Bra} (and is not the zero
space), there is an entire function $E$ (not unique) with the property:
\begin{equation}
  \label{eq:propE}
  \Imag(z)>0\implies |E(z)|>|E(\overline z)|
\end{equation} 
and in terms of which the evaluators are given as:
\begin{equation}\label{eq:eval}
  Z_z(w) = \frac{\overline{E(z)}E(w) - \overline{E^*(z)}E^*(w)}{i(\overline z-w)}
\end{equation}
We have written $E^*$ for the function $w\mapsto \overline{E(\overline w)}$.

Conversely to each  entire function satisfying \eqref{eq:propE} there is
associated a Hilbert space $H(E)$ with evaluators  given by \eqref{eq:eval}:
the elements of $H(E)$ are the entire functions $f$ such that both $f/E$ and
$f^*/E$ belong to the Hardy space of the upper half-plane ($f^*(z) = \overline{f(\overline z)}$).  A basic instance
of this theory is the Paley-Wiener space $H = PW_x$ ($x>0$) of entire
functions, square integrable on the real line, and  of exponential type at
most $x$. For these classical Paley-Wiener spaces we use the scalar product
$(f,f) = \frac1{2\pi}\int_\RR|f(z)|^2\,dz$, and the presence here of a
$\frac1{2\pi}$  is related to its absence in \eqref{eq:eval}; we find this
more convenient. An  appropriate $E$ function for $PW_x$ is $E(z) = e^{-ixz}$.

One of the axioms of \cite{Bra} is: (1) if $f\in H$ verifies $f(z_0) = 0$ then
$\frac{z-\overline{z_0}}{z - z_0} f(z)$ belongs to $H$ and has the same norm
as $f$. The others are: (2) $H$ is a Hilbert space and evaluating at a complex
point is a bounded linear form, and  (3) for any $f\in H$, the function
$f^*$
belongs to $H$ and has the same norm as $f$. From  axiom (1) one sees that
for $z_0$ non real, we can always find in $H$ (if it is not the zero space) a
function not vanishing at $z_0$. So the evaluator $Z_{z_0}$ can not be zero if
$z_0$ is non real (and if $H$ is not the zero space). Of course when $H$ is
known as an $H(E)$, then $\Imag(z)\neq 0\implies Z_z(z) = \|Z_z\|^2 >0$ by
\eqref{eq:eval} and \eqref{eq:propE}.

Let $\sigma = (z_1, \dots, z_n)$ be a finite sequence of non necessarily distinct complex
numbers with associated evaluators $Z_1$, \dots, $Z_n$ in $H$. More precisely,
in the case where the $z_i$'s  are not all distinct, we assemble in succession the indices corresponding to the same
complex number, and if for example $z_1 = z_2 = z_3 \neq z_4$, we let $Z_1 =
Z_{z_1}$, $Z_2 = Z_{z_1,1}$, $Z_3 = Z_{z_1,2}$, and $Z_4 = Z_{z_4}$ etc\dots
Also, when $f$ is an arbitrary analytic function, we introduce a notation
$f[z_j]$ such that in the example above $f[z_1] = f(z_1)$, $f[z_2] = f'(z_2)
(=f'(z_1))$, $f[z_3] = f''(z_3)$, $f[z_4] = f(z_4)$ etc\dots 

Let $H^\sigma$ be the closed subspace of $H$ of functions vanishing at the
$z_i$'s, in other words this is the orthogonal complement to the $Z_i$'s,
$1\leq i\leq n$.  For the classical Paley-Wiener spaces $H = PW_x$, such
subspaces $PW_x^\sigma$ have been considered in the work of Lyubarskii and
Seip \cite{LyuSeip}, where however the sequences $\sigma$ arising are
infinite. We restrict ourselves here to the finite case which is already
interesting: as will be shown in the companion paper \cite{painleve} this has
given a way from the classical  Paley-Wiener spaces to explicit Painlev\'e VI
transcendents. 

Let 
\begin{equation}
  \gamma(z) = \frac1{(z-z_1)\dots (z-z_n)}
\end{equation}
and define $H(\sigma) = \gamma(z) H^\sigma$:
\begin{equation}
  H(\sigma) = \{ F(z) = \gamma(z)f(z)\;|\; f\in H, f[z_1] = \dots = f[z_n] = 0\}
\end{equation} 

We call $F(z) = \gamma(z) f(z)$ the ``complete'' form of $f$ (for any $f$, in
$H$ or not, vanishing at the $z_i$'s), and call $z_1$, \dots, $z_n$ the
``trivial zeros''. We say that we switch from the space $H$ to the space
$H(\sigma)$ by adding trivial zeros, but this is of course slightly misleading
as the $z_i$'s are trivial zeros only for the incomplete functions $f(z)$, not
for the complete functions $F(z)$ which are vectors in the space $H(\sigma)$.

We give $H(\sigma)$ the Hilbert space structure which makes $f\mapsto F$ an
isometry with $H^\sigma$. Let us note that evaluations $F\mapsto F(z)$ are
again continuous linear forms on this new Hilbert space of entire functions:
this is immediate if $z\notin\sigma$ and follows from the Banach-Steinhaus
theorem if $z\in\sigma$. Let $F\in H(\sigma)$ with incomplete form $f$. If
$F(z_0) = 0$ then $f(z_0) = 0$, with multiplicity suitably increased  if $z_0$
belongs to $\sigma$. The function $g(z) = \frac{z-\overline{z_0}}{z - z_0}
f(z)$ belongs to $H$ and still vanishes on $\sigma$ (multiplicities included),
hence its complete form $\frac{z-\overline{z_0}}{z - z_0} F(z)$ belongs to
$H(\sigma)$. Finally, let us consider for $F\in H(\sigma)$ its conjugate in
the real axis $F^*(z) = \overline{F(\overline z)}$. With $F(z) = \gamma(z)
f(z)$ we thus have $F^*(z) = \gamma^*(z) f^*(z) = \gamma(z) \prod_{1\leq i
\leq n} \frac{z - z_i}{z - \overline{z_i}} f^*(z)$. But the function $f^*$
belongs to $H$ (with the same norm as $f$) and has  zeros at the
$\overline{z_i}$'s. Hence $\prod_{1\leq i \leq n} \frac{z - z_i}{z -
\overline{z_i}} f^*(z)$ belongs to $H$ with the same norm. And it has zeros at
the $z_i$'s, it is thus an element of $H^\sigma$ and its complete form is an
element of  $H(\sigma)$ (with the same norm as $F$). This completes the
verification that $H(\sigma)$ verifies the axioms of \cite{Bra} if $H$ does.

\begin{remark}
  Let us suppose that $z_0$ is not real. From what precedes if we can find a
  non-zero element $F$ in $H(\sigma)$ we can find one with $F(z_0)\neq0$. This
  proves in particular that if $Z_1$, \dots, $Z_n$ do not already span $H$,
  then any evaluator $Z_{z_0}$ with $z_0$ non-real and distinct from the
  $z_i$'s is not a linear combination of $Z_1$, \dots, $Z_n$.
\end{remark}

Let us give a first formula (which does not use \eqref{eq:eval} but only
\eqref{eq:ev}) for the evaluators $K_z$ in $H(\sigma)$ and their
scalar products $(K_w,K_z) = K_z(w)$. Let $k_z\in H^\sigma$ be the incomplete
form of $K_z$, so that $K_z(w) = \gamma(w) k_z(w)$. One has to be careful that
for $f\in H^\sigma$, with complete form $F$, we have by definition $(k_z,f) =
(K_z, F) = F(z) = \gamma(z)f(z) = \gamma(z) (Z_z,f)$. Hence $k_z$ is
$\overline{\gamma(z)}$ times the orthogonal projection $\pi(Z_z)$ of $Z_z\in
H$ onto $H^\sigma\subset H$. As is well-known, orthogonal projections can be
written in Gram determinantal form:
\begin{equation}
  k_z = \overline{\gamma(z)}\pi(Z_z) = \overline{\gamma(z)}\frac1{G_n}
  \begin{vmatrix}
    (Z_1,Z_1)&\dots &(Z_1,Z_n) &(Z_1,Z_z)\\
   (Z_2,Z_1)&\dots &(Z_2,Z_n) &(Z_2,Z_z)\\
\vdots&\dots&\vdots&\vdots\\
Z_1&\dots& Z_n&Z_z
  \end{vmatrix}
\end{equation}
We wrote $G_n$ for the principal $n\times n$ minor.
Then $(K_w,K_z) = K_z(w) = \gamma(w) k_z(w)$, and we have thus obtained,
writing now $K_z^\sigma$ for $K_z$:

\begin{proposition}\label{prop:1}
Let $H$ be a Hilbert space of entire functions with continuous evaluators
$Z_z$: $\forall f\in H\; f(z) = (Z_z,f)$. Let
$\sigma = (z_1,\dots,z_n)$ be a finite sequence of (non necessarily distinct) complex numbers
with associated evaluators $Z_1$, \dots, $Z_n$, assumed to be linearly
independent. Let $H(\sigma)$ be the Hilbert space of entire functions which are complete
forms of the elements of $H$ vanishing on $\sigma$. The evaluators of
$H(\sigma)$ are given by:
  \begin{equation}\label{eq:K}
    K_z^\sigma(w) = \frac{\gamma(w)\overline{\gamma(z)}}{G_n^\sigma}
    \begin{vmatrix}
      (Z_1,Z_1)&\dots &(Z_1,Z_n) &(Z_1,Z_z)\\
      (Z_2,Z_1)&\dots &(Z_2,Z_n) &(Z_2,Z_z)\\
      \vdots&\dots&\vdots&\vdots\\
      (Z_w,Z_1)&\dots& (Z_w,Z_n)&(Z_w,Z_z)
    \end{vmatrix}
  \end{equation}
where $G_n^\sigma$ is the principal $n\times n$ minor of the matrix at the
numerator. Of course this formula must be interpreted as a limit when $z$ or $w$
belongs to  $\sigma$.
\end{proposition}

Assume that   an $E$ function is known such that the evaluators
in $H$ are given by formula \eqref{eq:eval}.
We  find a function $E_\sigma$ playing the analogous role for $H(\sigma)$:
\begin{theorem}\label{thm:En}
  Let $E_\sigma(w)$ be the unique entire function such that its ``incomplete
  form'' (its product with $\prod_{1\leq i\leq
    n}(w-z_i)$) differs from $E(w)$ by a finite linear
  combination of the evaluators $Z_1(w)$, \dots, $Z_n(w)$. In other words, let
\begin{equation}\label{eq:E}
  E_\sigma(w) = \frac{\gamma(w)}{G_n^\sigma}
  \begin{vmatrix}
    (Z_1,Z_1)&\dots &(Z_1,Z_n) &E[z_1]\\
   (Z_2,Z_1)&\dots &(Z_2,Z_n) &E[z_2]\\
\vdots&\dots&\vdots&\vdots\\
Z_1(w)&\dots& Z_n(w)&E(w)
  \end{vmatrix}
\end{equation} where $G_n^\sigma$ is the principal $n\times n$ minor. The evaluator
$K_z^\sigma$ at $z$ for the space $H(\sigma)$ verifies:
\begin{equation}\label{eq:eval2}
  K_z^\sigma(w) = (K_w^\sigma,K_z^\sigma) = \frac{\overline{E_\sigma(z)}E_\sigma(w) - \overline{E_\sigma^*(z)}E_\sigma^*(w)}{i(\overline z-w)}\,
\end{equation}
\end{theorem}

\begin{remark}
  We mentioned earlier that, if $H$ is not already spanned by the $Z_i$,
  $1\leq i\leq n$, any evaluator $Z_z$ with $\Imag(z)\neq0$ is linearly
  independent from the $Z_i$'s. This implies $K_z^\sigma\neq 0$, hence for
  $\Imag(z)>0$ and by \eqref{eq:eval2}: $|E_\sigma(z)|^2 - |E_\sigma^*(z)|^2
  >0$, thus $E_\sigma$  given by \eqref{eq:E} indeed verifies \eqref{eq:propE} if $H(\sigma)$ is not
  the zero space.
\end{remark}

\begin{remark}
  Let us write $F(w) = E^*(w) = \overline{E(\overline w)}$ and similarly
  $F_\sigma(w) = E_\sigma^*(w) = \overline{E_\sigma(\overline w)}$. It will
  be shown in the proof that $F_\sigma$ follows the same recipe as $E_\sigma$:
  \begin{equation}
    \label{eq:F}
    F_\sigma(w) = \frac{\gamma(w)}{G_n^\sigma}
    \begin{vmatrix}
      (Z_1,Z_1)&\dots &(Z_1,Z_n) &F[z_1]\\
      (Z_2,Z_1)&\dots &(Z_2,Z_n) &F[z_2]\\
      \vdots&\dots&\vdots&\vdots\\
      Z_1(w)&\dots& Z_n(w)&F(w)
    \end{vmatrix}
  \end{equation} 
\end{remark}

\begin{remark} We did not see an immediate easy manipulation of determinants
leading to \eqref{eq:eval2} from \eqref{eq:K}  and \eqref{eq:E}. Even the
compatibility of the two equations \eqref{eq:E} and \eqref{eq:F} with the
relation $F_\sigma = E_\sigma^*$ does not seem to be immediately visible from
easy manipulations of determinants. However, under the additional hypotheses
that the space $H$ has the additional symmetry $f(z)\mapsto f(-z)$ and that
the ``trivial zeros'' are purely imaginary and distinct, a relatively simple
determinantal approach is proposed in \cite{painleve}. It leads in fact to
other determinantal expressions for $E_\sigma$ and $F_\sigma$ than
\eqref{eq:E} and \eqref{eq:F}, thus giving further determinantal
identities.
\end{remark}


\section{Adding one zero}

We establish the case $n=1$ of Theorem \ref{thm:En} by direct
computation. Appropriate notations are needed in order
to complete this deceptively simple looking task. We define:
\begin{align}
  \cE(w) &= \frac1{w-z_1}\left(E(w) - E(z_1)\frac{Z_1(w)}{(Z_1,Z_1)}\right)\\
  \cF(w) &= \frac1{w-z_1}\left(F(w) - F(z_1)\frac{Z_1(w)}{(Z_1,Z_1)}\right)\\
e_1 &= E(z_1)\qquad f_1 = F(z_1)
\end{align}
\begin{align}
\label{eq:a}  \text{hence: }\qquad  E(w) &= (w-z_1)\cE(w) + e_1 \frac{\overline{e_1}E(w) -
    \overline{f_1}F(w)}{i(Z_1,Z_1)(\overline{z_1} - w)}\\
\implies\; F(w) &= (w-\overline{z_1})\cE^*(w) - \overline{e_1} \frac{e_1 F(w) -
    {f_1}E(w)}{i(Z_1,Z_1)(z_1 - w)}\\
\label{eq:b}\text{on the other hand:}\qquad F(w) &= (w-z_1)\cF(w) + f_1 \frac{\overline{e_1}E(w) -
    \overline{f_1}F(w)}{i(Z_1,Z_1)(\overline{z_1} - w)}
\end{align}
We multiply the last identity by $\overline{z_1} - w$, the one before by $z_1 -w$
and substract:
\begin{equation}
  (z_1 - \overline{z_1})F(w) = (w-z_1)(\overline{z_1} - w)(\cE^*(w) - \cF(w)) -
  \frac{|e_1|^2 - |f_1|^2}{i(Z_1,Z_1)}F(w)
\end{equation}
Thus $(w-z_1)(\overline{z_1} - w)(\cE^*(w) - \cF(w)) = 0$ and we have established:
\begin{equation}\label{eq:EF}
  \cF(w) = \cE^*(w)
\end{equation}

We will also need the following identity:
\begin{equation}\label{eq:truc}
  \overline{e_1}\, \cE(w) - \overline{f_1}\, \cF(w) = -i Z_1(w)
\end{equation}
Indeed, from \eqref{eq:a} and \eqref{eq:b}:
\begin{equation}
  \overline{e_1} E(w) - \overline{f_1} F(w) =(w-z_1)(\overline{e_1} \cE(w) - \overline{f_1}
  \cF(w)) + (|e_1|^2 - |f_1|^2)\frac{Z_1(w)}{(Z_1,Z_1)}
\end{equation}
\begin{equation}
  \implies\; i(\overline{z_1}-w)Z_1(w) = (w-z_1)(\overline{e_1} \cE(w) - \overline{f_1}
  \cF(w))  + i(\overline{z_1}-z_1)Z_1(w)
\end{equation}
This proves \eqref{eq:truc}.

Let us now compute the determinant
\begin{equation}\label{eq:G}
      \begin{vmatrix}
      (Z_1,Z_1)&(Z_1,Z_z)\\
      (Z_w,Z_1)&(Z_w,Z_z)
    \end{vmatrix} =    
      (Z_1,Z_1)\frac{\overline{E(z)}E(w) - \overline{F(z)}F(w)}{i(\overline z -w)}
      - \overline{Z_1(z)}Z_1(w)
\end{equation}
We first consider:
\begin{align}
\label{eq:C}
   &  (Z_1,Z_1)\overline{E(z)}E(w) - (Z_1,Z_1)\overline{F(z)}F(w) - i(\overline{z_1}
    - w) \overline{Z_1(z)}Z_1(w)\\
    &= (Z_1,Z_1)\overline{E(z)}E(w) - (Z_1,Z_1)\overline{F(z)}F(w) -
    \overline{Z_1(z)}(\overline{e_1}E(w) - \overline{f_1}F(w)) \\
 &= (Z_1,Z_1)\overline{(z-z_1)}\left( \overline{\cE(z)} E(w) -  \overline{\cF(z)} F(w)\right) 
\end{align}
Using \eqref{eq:a} and \eqref{eq:b} this is equal to
\begin{align}
  &  (Z_1,Z_1){(\overline z-\overline{z_1})} \left( \overline{\cE(z)} (w-z_1)\cE(w) +
      \overline{\cE(z)}\frac{e_1 Z_1(w)}{(Z_1,Z_1)} - \overline{\cF(z)}
      (w-z_1)\cF(w) - \overline{\cF(z)}\frac{f_1 Z_1(w)}{(Z_1,Z_1)} \right)\\
&= (Z_1,Z_1)(\overline z-\overline{z_1})(w-z_1)( \overline{\cE(z)} \cE(w) - \overline{\cF(z)}
      \cF(w)) 
      + (e_1 \overline{\cE(z)} - f_1\overline{\cF(z)}) (\overline
      z-\overline{z_1})Z_1(w)\\
&= (Z_1,Z_1)(\overline z-\overline{z_1})(w-z_1)( \overline{\cE(z)} \cE(w) - \overline{\cF(z)}
      \cF(w)) 
      + i(\overline
      z-\overline{z_1})\overline{Z_1(z)}Z_1(w)\label{eq:D}
  \end{align}
where \eqref{eq:truc} was used. Identity of \eqref{eq:C} and \eqref{eq:D}
gives:
\begin{equation}
  \begin{split}
    (Z_1,Z_1)\overline{E(z)}E(w) - (Z_1,Z_1)\overline{F(z)}F(w) - i(\overline{z} -
    w) \overline{Z_1(z)}Z_1(w) \\= (Z_1,Z_1)(w-z_1)(\overline z - \overline{z_1})\left(
    \overline{\cE(z)} \cE(w) - \overline{\cF(z)} \cF(w)\right)
  \end{split}
\end{equation}
Comparison with \eqref{eq:G} gives the final result:
\begin{equation}
 \frac{\gamma_1(w)\overline{\gamma_1(z)}}{(Z_1,Z_1)} \begin{vmatrix}
      (Z_1,Z_1)&(Z_1,Z_z)\\
      (Z_w,Z_1)&(Z_w,Z_z)
    \end{vmatrix} = \frac{\overline{\cE(z)} \cE(w) - \overline{\cF(z)}
      \cF(w)}{i(\overline z -w)}
\end{equation}
As we know from \eqref{eq:EF} that $\cF = \cE^*$ this completes the proof of
Theorem \ref{thm:En} in the case $n=1$.

\section{General case with distinct trivial zeros}

An  induction establishes Theorem \ref{thm:En} when the $z_i$'s are
distinct. Let us suppose it true for the $n-1$ added ``trivial zeros'' $z_1$,
\dots ,$z_{n-1}$. From \eqref{eq:K} we know that for any
$z\in\CC\setminus\{z_1, \dots, z_{n-1}\}$, $\prod_{1\leq i\leq n-1}
(w-z_i)K_z^{z_1,\dots, z_{n-1}}(w)$ is a linear combination of the original
evaluators $Z_1(w)$, \dots, $Z_{n-1}(w)$ and $Z_z(w)$. This applies in
particular to $z = z_n$. The induction hypothesis tells us that $\prod_{1\leq
i<n} (w-z_i)E^{z_1,\dots, z_{n-1}}(w)$ differs from $E(w)$ by a linear
combination of the original evaluators $Z_1(w)$, \dots, $Z_{n-1}(w)$. The case
$n=1$ tells us that $(w-z_n)E^{z_1,\dots, z_{n}}(w)$ differs from
$E^{z_1,\dots, z_{n-1}}(w)$ by a multiple of $K_{z_n}^{z_1,\dots,
z_{n-1}}(w)$. Hence $\prod_{1\leq i\leq n} (w-z_i)E^{z_1,\dots, z_{n}}(w)$
differs from $E(w)$ by a linear combination of the original evaluators
$Z_1(w)$, \dots, $Z_{n-1}(w)$ and $Z_n(w)$, as was to be established. This
linear combination is fixed in a unique manner by evaluating at the trivial
zeros of $\prod_{1\leq i\leq n} (w-z_i)E^{z_1,\dots, z_{n}}(w)$. This
completes the proof of Theorem \ref{thm:En}. Furthermore we proved that the
same iterative recipe as for $E_\sigma$ works for the  construction of
$F_\sigma = E_\sigma^*$. Hence the formula \eqref{eq:F} holds.

\section{General case with multiplicities}

We introduce some notations for the case of repetitions among $z_1$, \dots,
$z_n$. We define $k_1$, \dots, $k_n$ such that $i' = i - k_i$ is the first
index with $z_{i'} = z_i$. For example if $z_1 = z_2 = z_3 \neq z_4 = z_5$,
$k_1 = 0$, $k_2 = 1$, $k_3 = 2$, $k_4 = 0$, $k_5 = 1$, etc\dots Then, for $f$
an analytic function, we recall the notation $f[z_i] := f^{(k_i)}(z_i)$. For
$f\in H$ we also have the scalar product $(Z_i,f) = f[z_i]$.

Let $k^\sigma(z,w)$ be the incomplete form of the
reproducing kernel $K^\sigma(z,w)$ in $H(\sigma)$:
\begin{equation}\label{eq:beta}
     k^\sigma(z,w) = \frac{1}{G_n}
    \begin{vmatrix}
      (Z_1,Z_1)&\dots &(Z_1,Z_n) &(Z_1,Z_z)\\
      (Z_2,Z_1)&\dots &(Z_2,Z_n) &(Z_2,Z_z)\\
      \vdots&\dots&\vdots&\vdots\\
      (Z_w,Z_1)&\dots& (Z_w,Z_n)&(Z_w,Z_z)
    \end{vmatrix} = Z_z(w)- \sum_{1\leq j \leq n} \beta_j^\sigma Z_j(w)
\end{equation} where the coefficients $\beta_1^\sigma$, \dots,
$\beta_n^\sigma$ (which are also functions of $z$) are determined by the constraints:
\begin{equation}
  \forall i\qquad \sum_{1\leq j \leq n} \beta_j^\sigma Z_j[z_i] =
  Z_z[z_i]\quad (=(Z_i,Z_z))
\end{equation}

Let 
$\epsilon>0$ and  $z_i^\epsilon = z_i - k_i \epsilon$ for $1\leq i\leq n$.
We let $k^\epsilon(z,w)$ be the (completely) incomplete form of the reproducing
kernel in $H(z_1^\epsilon,\dots, z_n^\epsilon)$. It thus has the shape:
\begin{equation}
  k^\epsilon(z,w) = Z_z(w) - \sum_{1\leq j \leq n} \alpha_j Z_{z_j^\epsilon}(w)
\end{equation}
with the constraints
\begin{equation}
  \forall i\qquad   k^\epsilon(z, z_i - k_i
  \epsilon) = 0
\end{equation}
Let us use the definitions
\begin{equation}
  Z_j^\epsilon := \epsilon^{-k_j} \sum_{0\leq m\leq k_j} (-1)^{m} {k_j \choose m}
  Z_{z_j - m \epsilon}
\end{equation}
to obtain vectors $Z_j^\epsilon \in H$ which span the same
   subspace of $H$ as  the $Z_{z_j^\epsilon} = Z_{z_j - k_j
  \epsilon}$. So there are coefficients $\beta_1^\epsilon$, \dots,
$\beta_n^\epsilon$ such that
\begin{equation}\label{eq:betaeps}
  k^\epsilon(z,w) = Z_z(w) - \sum_{1\leq j \leq n} \beta_j^\epsilon Z_j^\epsilon(w)
\end{equation}
with the constraints
\begin{equation}
  \forall i\qquad  \sum_{1\leq j \leq n} \beta_j^\epsilon
  Z_j^\epsilon(z_i - \epsilon k_i) = Z_z(z_i - \epsilon k_i)
\end{equation}
 We define the symbol, for any arbitrary analytic function on $\CC$:
\begin{equation}
  f[z_i^\epsilon] := \epsilon^{-k_i} \sum_{0\leq l\leq k_i} (-1)^{l} {k_i \choose l}
  f(z_i - l \epsilon) 
\end{equation}
With the help of these symbols, the constraints on the $\beta_j^\epsilon$ can
be equivalently rewritten:
\begin{equation}\label{eq:constraints}
  \forall i\qquad  \sum_{1\leq j \leq n} \beta_j^\epsilon
  Z_j^\epsilon[z_i^\epsilon] = Z_z[z_i^\epsilon]
\end{equation}
 We examine the behavior for $\epsilon\to0$ of the quantities
 $Z_j^\epsilon[z_i^\epsilon]$. In terms of the reproducing kernel $Z(z,w) =
 Z_z(w) = (Z_w,Z_z)$, which from equation \eqref{eq:eval}  is analytic in $w$ and
 anti-analytic in $z$ we have:
 \begin{equation}
   Z_j^\epsilon[z_i^\epsilon]  =
   \frac1{\epsilon^{k_i+k_j}}\sum_{
     \begin{smallmatrix}
       0\leq l\leq k_i\\0\leq m\leq k_j
     \end{smallmatrix}
}
(-1)^{l + m}{k_i \choose l}{k_j\choose m} Z(z_j - m\epsilon, z_i -
l\epsilon) 
\end{equation}
Thus, with $(\partial F)(w_1,w_2) = \frac{\partial}{\partial w_2} F(w_1,w_2)$, $(\delta
F)(w_1,w_2) = \frac{\partial}{\partial \overline{w_1}} F(w_1,w_2)$: 
\begin{equation}
  \lim_{\epsilon\to0} Z_j^\epsilon[z_i^\epsilon] = (\partial^{k_i}\delta^{k_j}
  Z)(z_j, z_i)
\end{equation}
On the other hand we have:
\begin{equation}
  \begin{split}
    (Z_i,Z_j) = \left.\frac{\partial^{k_i}}{\partial w^{k_i}}\right|_{w=z_i}
    Z_j(w) = \left.\frac{\partial^{k_i}}{\partial w^{k_i}}\right|_{w=z_i}
    \overline{(Z_j,Z_w)} = \left.\frac{\partial^{k_i}}{\partial
        w^{k_i}}\right|_{w=z_i} \overline{\left.\frac{\partial^{k_j}}{\partial
          \omega^{k_j}}\right|_{\omega=z_j} Z(w,\omega)}\\
= \left.\frac{\partial^{k_i}}{\partial
        w^{k_i}}\right|_{w=z_i} \left.\frac{\partial^{k_j}}{\partial
          \overline{\omega}^{k_j}}\right|_{\omega=z_j} Z(\omega,w) =  (\partial^{k_i} \delta^{k_j} Z)(z_j,z_i)
  \end{split}
\end{equation}
which gives
\begin{equation}
(Z_i,Z_j) =   (\partial^{k_i} \delta^{k_j} Z)(z_j,z_i) = \lim_{\epsilon\to0} Z_j^\epsilon[z_i^\epsilon] 
\end{equation}
There also holds, for any $z$:
\begin{equation}
  (Z_i,Z_z) = \lim_{\epsilon\to0} Z_z[z_i^\epsilon]
\end{equation}
So in the limit $\epsilon\to0$ the linear constraints \eqref{eq:constraints}
on $(\beta_j^\epsilon)_{1\leq j \leq n}$ become the constraints on the
coefficients $\beta_1^\sigma$,
\dots, $\beta_n^\sigma$ which give in \eqref{eq:beta} the incomplete reproducing kernel
$k^\sigma(z,w)$. This shows in passing that for $\epsilon\neq0$ small the
vectors $Z_j^\epsilon$ are also linearly independent, and proves $\forall j\;
\lim_{\epsilon\to0} \beta_j^\epsilon = \beta_j^\sigma$. We have further
\begin{equation}
Z_j(w) = \lim_{\epsilon\to0} Z_j^\epsilon(w)
\end{equation}
and this finally establishes:
\begin{equation}
  \lim_{\epsilon\to0} k^\epsilon(z,w) = k^\sigma(z,w)
\end{equation}

In the exact same manner we can examine the quantities:
\begin{equation}\label{eq:Eincomplete}
  e_\sigma(w) = \frac{1}{G_n^\sigma}
  \begin{vmatrix}
    (Z_1,Z_1)&\dots &(Z_1,Z_n) &E[z_1]\\
   (Z_2,Z_1)&\dots &(Z_2,Z_n) &E[z_2]\\
\vdots&\dots&\vdots&\vdots\\
Z_1(w)&\dots& Z_n(w)&E(w)
  \end{vmatrix}
\end{equation}
and
\begin{equation}\label{eq:Eincomplete}
  e_\epsilon(w) = \frac{1}{G_n^\epsilon}
  \begin{vmatrix}
    (Z_{z_1^\epsilon},Z_{z_1^\epsilon})&\dots &(Z_{z_1^\epsilon},Z_{z_n^\epsilon}) &E({z_1^\epsilon})\\
   (Z_{z_2^\epsilon},Z_{z_1^\epsilon})&\dots &(Z_{z_2^\epsilon},Z_{z_n^\epsilon}) &E({z_2^\epsilon})\\
\vdots&\dots&\vdots&\vdots\\
Z_{z_1^\epsilon}(w)&\dots& Z_{z_n^\epsilon}(w)&E(w)
  \end{vmatrix}
\end{equation}
and prove 
\begin{equation}
  \lim_{\epsilon\to0} e_\epsilon(w) = e_\sigma(w)
\end{equation}
There is also an immediate limit to be taken in the gamma factor, and in the
end we obtain the reproducing kernel formula \eqref{eq:eval2} for the space
$H(\sigma)$. The formula $E_\sigma^* = F_\sigma$ with $F_\sigma$ given by
$\eqref{eq:F}$ is proven in the same manner.

\section{An example}

We take $H = PW_x$ ($x>0$), the Paley-Wiener space of
entire functions of exponential type at most $x$, square integrable on the
real line ($(f,f) = \frac1{2\pi}\int_\RR|f(z)|^2\,dz$). The evaluators are
given by the formula
\begin{equation}
  Z_z(w) = (Z_w,Z_z) = 2\frac{\sin((\overline z-w)x)}{\overline z -w} \quad (=\int_{-x}^x e^{iwt}e^{-i\overline zt}\,dt)
\end{equation}
We can choose $E(z) = e^{-ixz}$, $F(z) = E^*(z) = e^{+ixz}$. Let $z_1$, $z_2$,
\dots, $z_n$ and $z$ be distinct complex numbers and define

\begin{equation}
  G_{n} =     \begin{vmatrix}
     \strut  2\frac{\sin((\overline{z_1}-z_1)x)}{\overline{z_1}-z_1}&2\frac{\sin((\overline{z_2}-z_1)x)}{\overline z_2 -z_1}&\dots&2\frac{\sin((\overline{z_n}-z_1)x)}{\overline{z_n} -z_1}\\
     \mathstrut  2\frac{\sin((\overline{z_1}-z_2)x)}{\overline z_1 -z_2}&2\frac{\sin((\overline{z_2}-z_2)x)}{\overline{z_2}-z_2}&\dots &2\frac{\sin((\overline{z_n}-z_2)x)}{\overline{z_n} -z_2} \\
     \mathstrut  \vdots&\dots&\dots&\vdots\\
     \mathstrut  2\frac{\sin((\overline{z_1}-z_n)x)}{\overline{z_1} -z_n}&\dots&\dots  &2\frac{\sin((\overline{z_n}-z_n)x)}{\overline{z_n}-z_n}
    \end{vmatrix}
\end{equation}

\begin{equation}
  G_{n}(z,z) =     \begin{vmatrix}
      2\frac{\sin((\overline{z_1}-z_1)x)}{\overline{z_1}-z_1}&2\frac{\sin((\overline{z_2}-z_1)x)}{\overline{z_2} -z_1}&\dots &2\frac{\sin((\overline{z_n}-z_1)x)}{\overline{z_n} -z_1} &2\frac{\sin((\overline{z}-z_1)x)}{\overline{z} -z_1}\\
      2\frac{\sin((\overline{z_1}-z_2)x)}{\overline{z_1} -z_2}&2\frac{\sin((\overline{z_2}-z_2)x)}{\overline{z_2}-z_2}&\dots &2\frac{\sin((\overline{z_n}-z_2)x)}{\overline{z_n} -z_2} &2\frac{\sin((\overline{z}-z_2)x)}{\overline{z} -z_2}\\
      \vdots&\dots&\vdots&\vdots&\vdots\\
      2\frac{\sin((\overline{z_1}-z)x)}{\overline{z_1} -z}&\dots&\dots  &2\frac{\sin((\overline{z_n}-z)x)}{\overline{z_n}-z} &2\frac{\sin((\overline{z}-z)x)}{\overline{z}-z}
    \end{vmatrix}
\end{equation}

\begin{equation}
  e_{n}(z) =     \begin{vmatrix}
       2\frac{\sin((\overline{z_1}-z_1)x)}{\overline{z_1}-z_1}&2\frac{\sin((\overline{z_2}-z_1)x)}{\overline{z_2} -z_1}&\dots &2\frac{\sin((\overline{z_n}-z_1)x)}{\overline{z_n} -z_1} &e^{-ixz_1}\\
      2\frac{\sin((\overline{z_1}-z_2)x)}{\overline{z_1} -z_2}&2\frac{\sin((\overline{z_2}-z_2)x)}{\overline{z_2}-z_2}&\dots
      &2\frac{\sin((\overline{z_n}-z_2)x)}{\overline{z_n} -z_2} &e^{-ix z_2}\\
      \vdots&\dots&\vdots&\vdots&\vdots\\
      2\frac{\sin((\overline{z_1}-z)x)}{\overline{z_1} -z}&\dots&\dots  &2\frac{\sin((\overline{z_n}-z)x)}{\overline{z_n}-z} &e^{-ixz}
    \end{vmatrix}
\end{equation}

\begin{equation}
  f_{n}(z) =     \begin{vmatrix}
       2\frac{\sin((\overline{z_1}-z_1)x)}{\overline{z_1}-z_1}&2\frac{\sin((\overline{z_2}-z_1)x)}{\overline{z_2} -z_1}&\dots &2\frac{\sin((\overline{z_n}-z_1)x)}{\overline{z_n} -z_1} &e^{ixz_1}\\
      2\frac{\sin((\overline{z_1}-z_2)x)}{\overline{z_1} -z_2}&2\frac{\sin((\overline{z_2}-z_2)x)}{\overline{z_2}-z_2}&\dots
      &2\frac{\sin((\overline{z_n}-z_2)x)}{\overline{z_n} -z_2} &e^{ix z_2}\\
      \vdots&\dots&\vdots&\vdots&\vdots\\
      2\frac{\sin((\overline{z_1}-z)x)}{\overline{z_1} -z}&\dots&\dots  &2\frac{\sin((\overline{z_n}-z)x)}{\overline{z_n}-z} &e^{ixz}
    \end{vmatrix}
\end{equation}

Then 
\begin{equation}
  {G_n(z,z)}{G_n} = \frac{|e_n(z)|^2 - |f_n(z)|^2}{2\,\text{Im}(z)}
\qquad\text{and}\qquad
  f_n(z) = \prod_{1\leq i\leq n} \frac{z-z_i}{z - \overline{z_i}}\;
  \overline{e_n(\overline z)}
\end{equation}


\end{document}